\documentclass[%
 reprint,
 amsmath,amssymb,
 aps,
pre,
]{revtex4-2}

\usepackage{graphicx}
\usepackage{dcolumn}
\usepackage{bm}
\usepackage{subcaption}

\begin{document}

\preprint{APS/123-QED}

\title{Do triangles matter? Replicating hypergraph disease dynamics with lower-order interactions}

\author{Eugene Tan}
    \altaffiliation[Also at ]{The Complex Systems Group, The University of Western Australia.}
    \affiliation{The Kids Research Institute Australia, Nedlands 6009, Australia}
    \email{eugene.tan@thekids.org.au}
\author{Michael Small}
    \email{michael.small@uwa.edu.au}
\author{Shannon D. Algar}%
    \email{shannon.algar@uwa.edu.au}
    \affiliation{The Complex Systems Group, The University of Western Australia, Crawley 6009, Australia}

\date{\today}

\begin{abstract}
Disease spreading models such as the ubiquitous SIS compartmental model and its numerous variants are widely used to understand and predict the behaviour of a given epidemic or information diffusion process. A common approach to imbue more realism to the spreading process is to constrain simulations to a network structure, where connected nodes update their disease state based on pairwise interactions along the edges of their local neighbourhood. Simplicial contagion models (SCM) extend this to hypergraphs such that groups of three nodes are able to interact and propagate the disease along higher-order hyperedges (triangles). Though more flexible, it is not clear the extent to which the inclusion of these higher-order interactions result in dynamics that are characteristically different to those attained from simpler pairwise interactions. Here, we propose an agent-based model that unifies the classical SIS/SIR compartmental model and SCM, and extends it to allow for interactions along hyperedges of arbitrary order. Using this model, we demonstrate how the steady-state dynamics of pairwise interactions can be made to replicate those of simulations that include higher-order topologies by linearly scaling disease parameters based on a proposed measure of network activity. By allowing disease parameters to dynamically vary over time, lower-order pairwise interactions can be made to closely replicate both the transient and steady-state dynamics of higher-order simulations. We demonstrate that this relationship is robust to misspecification in the assumed higher-order interaction model, and applies to non-clique complex hypergraphs with non-trivial heterogeneous topology. For the latter case, it is found that heterogeneities in hypergraph topology result in weakened approximations of higher-order dynamics by pairwise interactions.
\end{abstract}

\maketitle


\section{\label{sec:level1_intro}Introduction}

Spreading processes over a given topology are ubiquitous across many social and physical systems. From epidemics to information diffusion, many studies aim to understand the governing rules of spreading processes in social contexts \cite{pastor2015epidemic, jalili2017information, ferraz2024contagion, daley1964epidemics}. Nevertheless, the influence of the underlying topology on which spreading occurs can sometimes be as important as the mechanistic rules governing the spreading process \cite{iacopini2019simplicial,shirley2005impacts}. This consideration is crucial in the study of contagion and epidemics on social networks, where local properties like degree and clustering can mediate non-trivial dynamics such as super-spreader events and endemic states \cite{small2006super,shirley2005impacts}. 

Recent approaches propose the inclusion of interactions along higher-order topologies such as simplices (e.g. triangles, tetrahedrons) to account for social interactions that include more individuals than pairwise interactions \cite{iacopini2019simplicial, battiston2020networks, battiston2021physics}. The inclusion of higher-order interactions have been shown to promote critical transitions and bistability in epidemics \cite{zhao2024susceptible, lucas2023simplicially, bodo2016sis, ferraz2024contagion}. However, the importance that they play in altering the overall epidemic trajectory remains an active area of research. Current work has also explored the effect of simplicial contagion processes on infection paths on hypergraphs \cite{contreras2024infection}. Several analytical methods have also been developed to better understand the governing factors of unique simplicial contagion behaviours such as hysteresis, bistability and reduced epidemic thresholds \cite{kiss2025decoding, zhao2025higher}

In this work, we present a discrete numerical model that generalises the classical SIS/SIR compartmental model with homogeneous mixing to account for hyperedge interactions of arbitrary order. Comparing dynamics with and without higher-order interactions, we find topology primarily affects epidemic transients. If the network structure is known, one can calculate normalisation values related to the potential network activity. These values may be used to define disease parameters for first-order pairwise interactions that replicate the steady state dynamics of those that include higher-order interactions, but are insufficient for reproducing epidemic transients. However, we find a duality where the effects of network topology can be overcome by allowing for temporally varying disease parameters; allowing first-order dynamics to closely approximate both transient and steady state trajectories of higher-order systems. 

In addition to requiring knowledge of the order scaling of interaction dynamics, the normalisation approach makes two assumptions regarding the hypergraph topology. Namely, that (1) all hypergraphs are defined as clique complexes, (2) macro-scale topological features are uniformly distributed across the network. To investigate the robustness of the results, tests are conducted on cases where the scaling function is incorrect. We show robustness to model misspecification provided the size scales of the maximal infection amplification between the true and misspecified scaling function are of similar order. These results are broadly aligned with the recent analytical results by \cite{malizia2025pair} where mean field approximations of pairwise interactions were able to captured features of simplicial contagion. We also extend analyses to artificial and empirical networks that do not adhere to Assumptions (1) and (2) and show robustness of results when disease parameters are allowed to temporally vary. Further inspection of the network structure and component hyperedge activities suggest that non-trivial macro-scale heterogeneities in hypergraph topologies limit the approximation of higher-order dynamics by pairwise interactions.

\section{\label{sec:level1_model}Disease Model}
We begin by presenting a variation of the classical compartmental ODE model \cite{kermack1927contribution} focusing on the three compartment case with susceptible-infected-removed (SIR) states,
\begin{subequations}
\begin{align}
    \frac{\mathrm{d}S}{\mathrm{d}t} &= -\beta SI + \mu I\\
    \frac{\mathrm{d}I}{\mathrm{d}t} &= \beta SI - \mu  I - \alpha I \\
    \frac{\mathrm{d}R}{\mathrm{d}t} &= \alpha I
\end{align}
\end{subequations}
where $N = S+I+R$ is the total population and $S,I,R$ are the susceptible, infected and removed populations respectively. The constants $\beta$, $\mu$ and $\alpha$ are disease parameters corresponding to the rates of infection, recovery and removal. For $\alpha=0$, the SIR system reduces to the simpler SIS system. 

The simple SIR model assumes that populations homogeneously mix. However, this is not true for human populations where disease spread occurs on a social contact network. Here, individuals -- represented as nodes -- interact within a small neighbourhood with occasional branching, clustering or long range connections \cite{boccaletti2023structure, iacopini2019simplicial, kiss2017mathematics}. One can simulate more realistic spread dynamics by employing an agent-based approach where randomly selected nodes interact with all or a collection of their neighbours every time step \cite{pastor2001dynamical}. Whilst this node-based approach is more realistic for modelling the diffusion of abstract quantities such as information, it is not necessarily adequate for disease contagion as it treats individuals rather than interactions as facilitators of infection. It is unlikely for an individual to interact with all of their neighbours in a given time step. Instead, we propose simulating diseases on networks using an edge-based approach. 

Let $G=(E,V)$ be a hypergraph consisting of a set of hyperedges $E$ and nodes $V$. A $k^{th}$ order hyperedge is a connection between $k+1$ different nodes. Each simulation time step consists of $m$ iterations. A random hyperedge $e\in E$ is selected in each iteration. Each susceptible node $s_{i}$ in $e$ has a probability $\beta_e$ of getting infected and transitioning to state $I$,
\begin{equation}
    \beta_{e} = \min\left[f\left(\frac{n_{I}(e)}{n_{S}(e)-1} \right)n_{I}(e)\beta_{ABM},1 \right],
\end{equation}
where $n_I(e)$ and $n_{S}(e)$ are the number of nodes in $e$ that are infected and susceptible respectively, $\beta_{ABM}$ is a base infection probability for a pairwise interaction, and $f$ is a scaling function that controls the extent to which multiple infected nodes in a given hyperedge amplify the probability of infection. In simplicial social contagion models (SCM) where interactions occur up to order $k=2$, susceptible nodes whose neighbours are all infected experience an amplified probability of infection \cite{iacopini2019simplicial}. Whilst SCM has only been applied to social dynamics (e.g. knowledge spreading), we argue that a similar principal is also relevant when studying diseases where the infection probability may be mediated by the concentration of disease particles within an individual's immediate vicinity during an interaction. The allowance for variable infectivity with has been explored prior by Fu et al. \cite{fu2008epidemic}. Inspired by this, we use the following scaling function to describe the disease spread dynamics on hyperedges,
\begin{equation}
    f(p) = 
    \begin{cases}
        2p, &0\leq p<0.5 \\
        1+2k_{\gamma}(2p-1), & p \geq 0.5 
    \end{cases}
\end{equation}
where $k_{\gamma}$ is the maximum amount that a hyperedge interaction can amplify the infection probability. After interaction, all infected nodes recover or die with probabilities $P(I\to S) = \mu_{ABM}$ and $P(I\to R) = \alpha_{ABM}$.

\begin{figure}
    \centering
    \includegraphics[width=\linewidth]{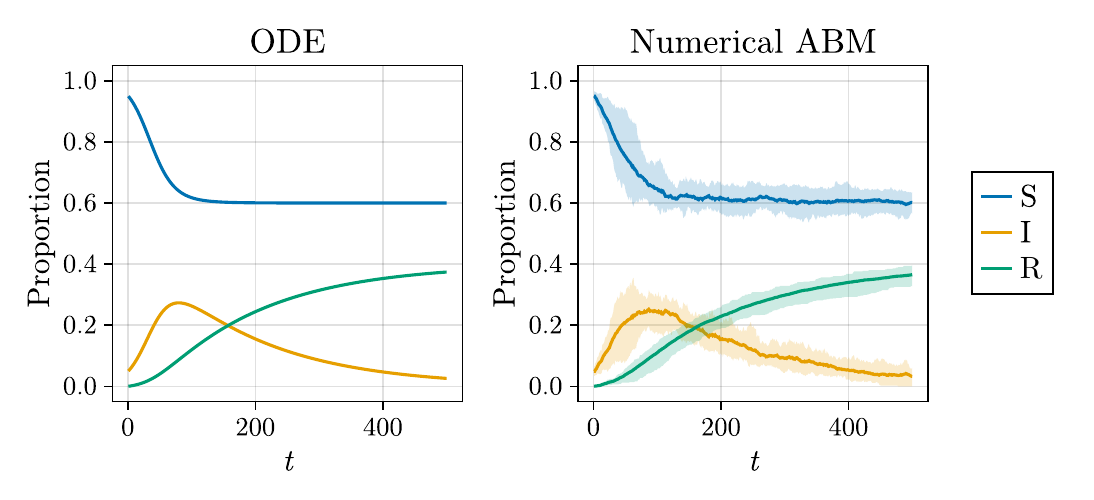}
    \caption{Epidemic trajectories from ODE model assuming homogeneous mixing, and 20 simulations of the numerical approximation using a $N=500$ node fully connected network. Mean and 90\% quantiles shown. Parameters are $\beta_{ODE}=0.0003$, $\mu_{ODE}=0.5N\beta_{ODE}$, $\alpha_{ODE} = 0.04N\beta_{ODE}$.}
    \label{fig:analytic_abm}
\end{figure}

The proposed numerical model achieves two goals. Firstly, it unifies the classical homogeneous ODE disease model with existing SCM up to order $k_{max}=2$. For $k_{max}=1$ and a fully connected graph $G$, this model is a numerical approximation for the SIR/SIS model using the following disease parameters,
\begin{subequations}
    \begin{align}
        \beta_{ABM} &= \frac{N(N-1)}{2} \beta_{ODE}\,\delta t, \\
        \mu_{ABM} &= \mu_{ODE}\, \delta t, \\
        \alpha_{ABM} &= \alpha_{ODE}\,\delta t,
    \end{align}
\end{subequations}
where $\delta t = 1/m$ is the integration time step. We validate simulations against the analytical ODE model and find good agreement in epidemic trajectories (see Fig. \ref{fig:analytic_abm}). For $k_{max} = 2$, the model is equivalent to an edge-based formulation of the simplicial contagion model (SCM) by Iacopini et al. \cite{iacopini2019simplicial}, with the conversion,
\begin{equation}
    k_{\gamma}= \frac{\beta +\beta_{\Delta}}{2\beta},
\end{equation}
where $\beta_{\Delta}$ is the simplicial infection rate. Secondly, this generalises simplicial contagion to account for simplices of arbitrary order and varying proportions of infected individuals within each hyperedge.

\section{\label{sec:level1_activities}Network Activities}
Identifying dynamical differences between networks with and without interactions on higher-order topologies requires an appropriate normalisation that allows for comparison. Thus, we consider two candidate quantities related to the rate of new infection $\frac{dI}{dt}$.

The first quantity, called the combinatorial network activity $\hat{\lambda}^{(K)}$ calculated up to order $K$ for a hypergraph $G$ with respect to a baseline infection rate $\beta$ is defined as,
\begin{equation}
    \hat{\lambda}^{(K)} (G,\beta)=\frac{\sum_{k=1}^{K} N_{k} \hat{\beta}^{(k)} }{\sum_{k=1}^{K} N_{k}},
\end{equation}
where
\begin{equation}
    \hat{\beta}^{(k)} = \frac{\sum_{i=0}^{k+1} \binom{k+1}{i} (k+1-i) \min\left[f\left(\frac{i}{k} \right)\,i\,\beta,1 \right]}{\sum_{i=0}^{k+1} \binom{k+1}{i}},
\end{equation}
and approximates $\frac{dI}{dt}$ under the assumption that individuals are equally likely to be infected or susceptible. Therefore, $\hat{\lambda}^{(K)}$ depends only on network the network topology and the scaling function $f$.

Similarly, an exact network activity $\lambda^{(K)}$ up to order $K$ can be calculated as,
\begin{equation}
    \lambda^{(K)}(G,\beta,X) = \frac{\sum_{e\in E_{K}} n_{S}(e) \beta_{e}}{|E_{K}|},
\end{equation}
where $E^{(K)}$ is the set of all hyperedges of order up to $k_{max}=K$, and $X$ are the nodes' states. Unlike $\hat{\lambda}^{(K)}$, $\lambda^{(K)}$ does not assume a distribution of node states and is instead a function of the full network state and approximates the instantaneous value of $\frac{dI}{dt}$.

\section{\label{sec:level1_results}Results}
\subsection{Initial Network Activity Normalisation}\label{subsec:init_norm}
To identify the role of higher-order interactions in disease spreading dynamics, we pose the following question: Can pairwise interactions replicate higher-order dynamics using an appropriate normalisation?

We address this question using the case of SIS dynamics with baseline disease parameters $(\beta,\mu,\alpha) = (0.05,0.0001,0)$. The SIS case is analysed first as it possesses steady states on which measures can be defined. We later extend our analyses to the more general SIR case. 

Given an undirected $N$ node network with a pairwise adjacency matrix $A$, we define two different hypergraphs $G_{1} = (E_{1},V)$, and $G_{K}=(E_{K},V)$ with hyperedges defined up to the order 1 and $K$ clique complex of $A$ respectively. Parallel simulations with the same initial conditions are conducted using disease parameters $(\beta_{1},\mu,\alpha)$ and $(\beta_{K},\mu,\alpha)$ with infection values normalised using the ratio of both networks' combinatorial network activity,
\begin{subequations}
    \begin{align}
        \beta_{1} & = \frac{\hat{\lambda}^{(K)}(G_{K},\beta)}{\hat{\lambda}^{(1)}(G_{1},\beta)} \beta_{K}, \\
        \beta_{K} & = \beta.
    \end{align}
\end{subequations}
This normalisation ensures that simulations on both $G_{1}$ and $G_{K}$ begin with approximately similar network activity. Thus, deviations in the epidemic trajectories of both simulations are the result of topology rather than disease parameters.

\begin{figure}
    \centering
    \includegraphics[width=0.95\linewidth]{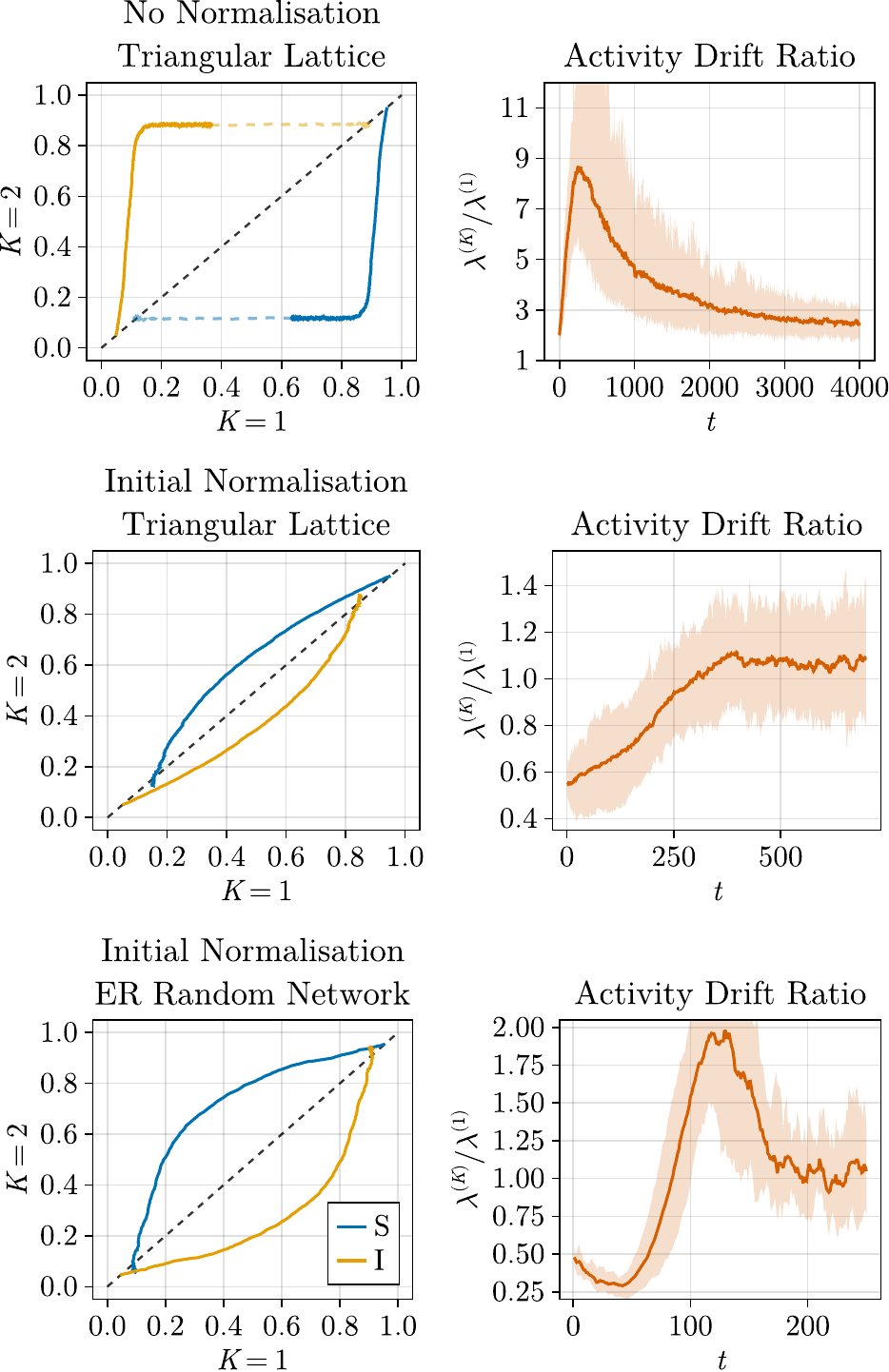}
    \caption{SI epidemic trajectories between the pairwise and higher-order simulations. Accompanying activity ratios are shown right revealing a saturation to 1 as the epidemic reaches the endemic steady state. Dotted trajectory in the unnormalised case corresponds to enabling of higher-order interactions in pairwise simulation.}
    \label{fig:init_norm}
\end{figure}

This test is run on two $N=500$ node networks with different topologies: a 2D triangular lattice, and an Erd\"{o}s-R\'{e}nyi (ER) random network with connection probability $p = 0.1$. The triangular lattice represents the case where the probability of higher order interactions is greatest (i.e. $\frac{|E_{K}|}{|E_{1}|}$ is maximised) for $K=2$. The random network is the case that allows for multiple orders of interaction, with maximum order chosen to be $K=4$. We run $50$ and $10$  simulations of the triangular lattice and random network respectively with a resolution of $m=50$ (see Fig. \ref{fig:init_norm}). 

A base case with unnormalised infection rates is run as a control. Enabling higher-order interactions after the pairwise system settles pushes the steady state align with the higher-order scenario. For the test case, applying a normalisation based on the ratio of each network's combinatorial network activity results in epidemic trajectories whose steady state in the $K=1$ case closely mimic those of $K>1$. However, the initial adjustment of disease parameters is insufficient for replicating epidemic transients in earlier portions of the trajectory. These findings suggest that network topology and higher-order interactions are critical during epidemic transients, but washes out in the endemic steady state. In early stages of the epidemic, infection spread is primarily mediated by the local network structures in the neighbourhood of infected agents, which in turn govern the slopes of the trajectory. Once steady state is achieved, pairwise first-order $K=1$ interactions given an appropriate normalisation is sufficient to reproduce spread dynamics. 

\subsection{Adaptive Network Activity Normalisation}\label{subsec:adaptive_norm}
The capacity for networks with the same initial activity to reproduce epidemic steady states naturally leads the discussion to a second related question: Under what minimal conditions can pairwise interactions be used to replicate the \textit{full dynamics} as those of higher-order? We again perform similar tests consisting of two parallel simulations with separate disease parameters with the exception that the infection rate of the $K=1$ case is allowed to vary over time, with parameters $(\beta_{1}(t), \mu, \alpha)$ and $(\beta_{K}, \mu, \alpha)$. This is plausible in cases where disease strains evolve, or intervention policies are enforced (i.e. dynamic network structures, altered infection parameters).

For simulating adaptive disease parameters, we first calculate the moving average of the exact network activity ratios in both hypergraphs,
\begin{subequations}
    \begin{align}
        \xi(t)  & = \frac{\lambda^{(K)}(G_{K},\beta,X_{K}(t))}{\lambda^{(1)}(G_{1},\beta_{1}(t),X_{1}(t))},\\
        \bar{\xi}(t)  & = \frac{1}{\tau}\sum_{n=t}^{t+\tau} \xi(n),
    \end{align}
\end{subequations}

where $\tau=5$ a set window length. The pairwise infection rate parameter $\beta_{1}(t)$ is dynamically adjusted with respect to a conservative set threshold $\rho=0.5$
\begin{equation}
    \beta_{1}(t) = 
    \begin{cases}
        \beta_{1}(t-1), & |\bar{\xi}(t)-1| < \rho \\
        \xi(t)\beta_{1}(t-1),  &  |\bar{\xi}(t)-1| \geq \rho
    \end{cases}
\end{equation}
The adjustment of $\beta_{1}(t)$ based on a threshold of the moving average $\bar{\xi}(t)$ limits the extent to which the infection parameter is allowed to vary, and thus aims to mimic a more realistic evolution of disease dynamics where base infection rates are less volatile.  

\begin{figure*}
    \centering
    \includegraphics[width=0.9\linewidth]{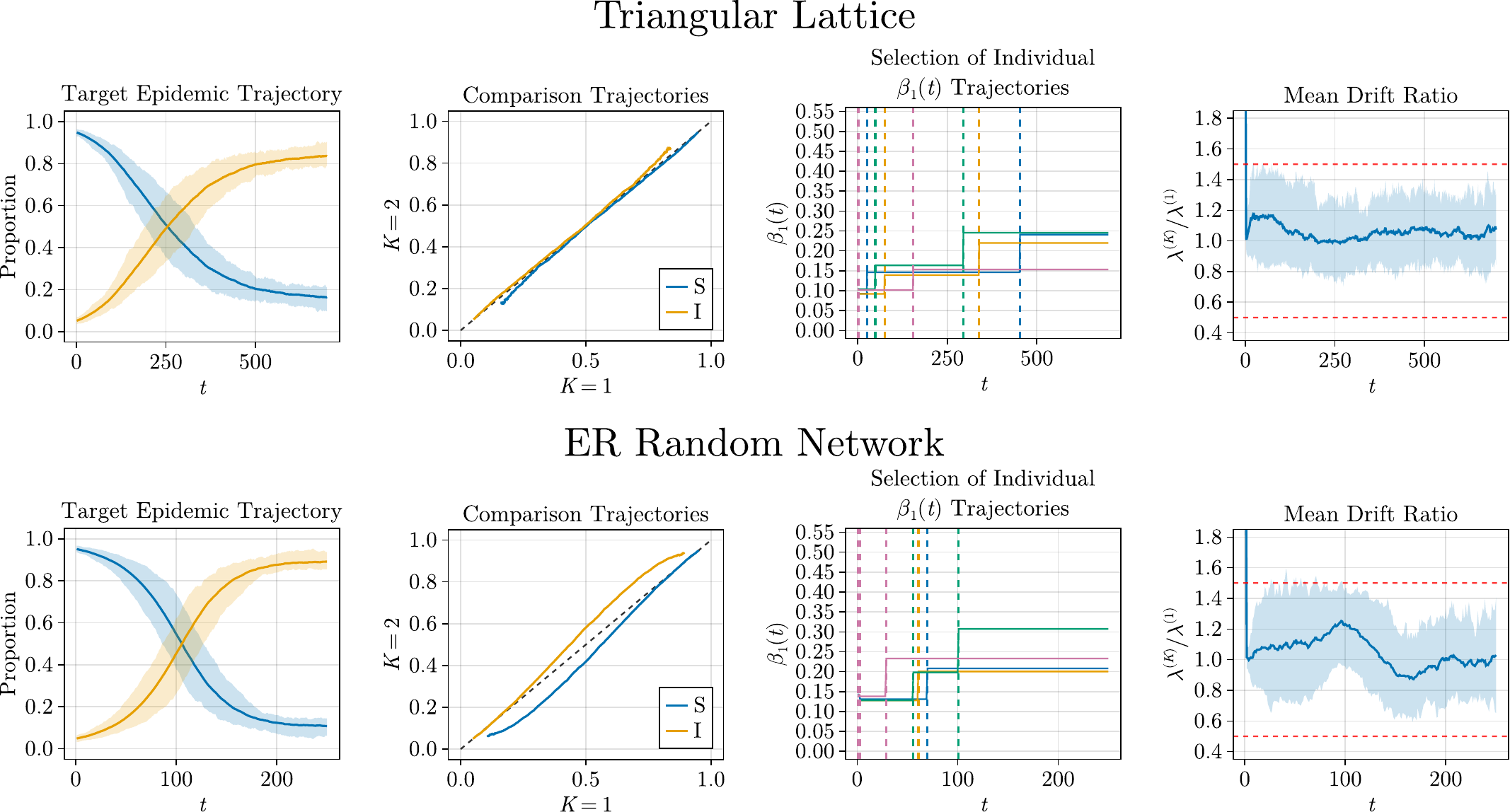}
    \caption{Simulations with dynamic disease parameters with $K=1$ for triangular lattice and random networks. Shaded regions correspond to $90\%$ CI. Left to right: target trajectories from the higher order simulation, comparison trajectories between pairwise and higher order case, four randomly chosen trajectories of $\beta_1(t)$ with adjustment events indicated by vertical lines, and the mean drift ratio $\bar{\xi}(t)$, with the $\rho$ threshold envelope in red.}
    \label{fig:adaptive}
\end{figure*}

Results for the normalisation tests with dynamically varying disease parameters are shown in Fig. \ref{fig:adaptive}. We find that the dynamical adjustment of disease parameters allows for first-order (pairwise) simulations to replicate those of higher order in both transient and steady state regimes. Furthermore, we note that the adjustments of disease parameters are infrequent, with 2-4 adjustment events occurring across 700 simulated time steps in each simulation. A majority of adjustments occur during the first half of the transient regime. These results suggest a duality between network topology and disease parameters. Namely, the effects of higher-order interactions may be partially accounted for in first-order approximations through the presence of dynamic disease parameters. 

\subsection{SIR Dynamics}
All of the previous tests focus on the SIS disease scenario, which is characterised by an epidemic transient followed by a steady state regime where the disease is either eradicated or remains in endemic. A more interesting, and arguably challenging case, is that of SIR where infection trajectories exhibit a single peak during the epidemic's transient phase. 

\begin{figure*}
    \centering
    \includegraphics[width=\linewidth]{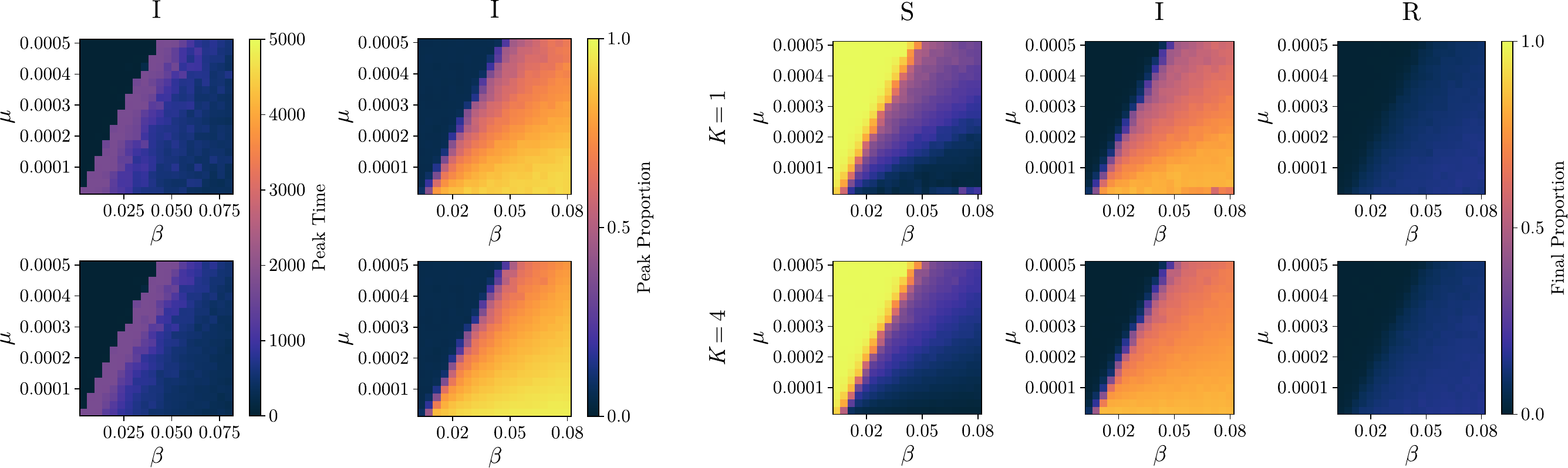}
    \caption{Phaseplots of the peak time, and peak and final proportions for simulated SIR model with dynamic disease parameters. Top and bottoms rows correspond to the pairwise $K=1$ and higher-order $K=4$ simulations for $N=500$ node ER random network. Dynamic adjustment shows good replication of phaseplots and epidemic outbreak threshold. }
    \label{fig:phaseplot}
\end{figure*}

We extend our analyses to this case and test the dynamic adjustment of infection rates applied across a range of base disease parameter values $\beta_{K} \in [0.004,0.08]$ and $\mu_{K} \in [0.000025,0.0005]$ with maximum order $K=4$ and $k_{\gamma}=3$. Phaseplots of the peak time, and peak and final proportions are given in Fig. \ref{fig:phaseplot}. Overall, we find that normalisations using activity ratios enable higher-order interaction dynamics to be replicated by those of lower order.

\subsection{Robustness and Model Misspecification}
\begin{figure}
    \centering
    \includegraphics[width=0.6\linewidth]{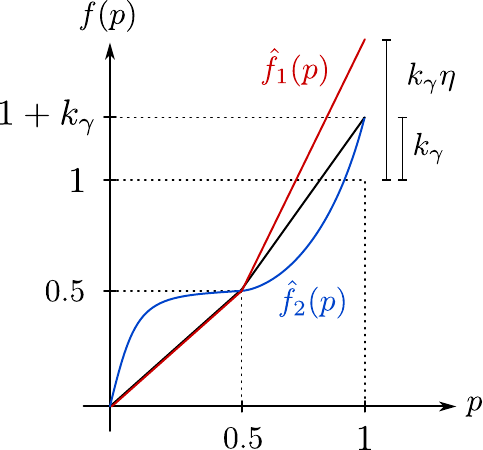}
    \caption{Schematic of scaling function $f(p)$ and tested misspecified forms $\hat{f}_1(p)$ and $\hat{f}_2(p)$.}
    \label{fig:scaling_functions}
\end{figure}

While the thresholds imposed for the adaptive adjustment of disease parameters are relatively mild, the network activity normalisation method requires extensive knowledge about the scaling of spreading dynamics up to the $K^{th}$ order hyperedge. This is not easily defined in practice, and direct application will almost certainly contain errors in the specification of the scaling function $f(p)$. To test the robustness of our findings, the same numerical experiment is repeated using a misspecified scaling function $\hat{f}(p)$ when calculating activity ratios $\xi$. We test two types of misspecification: scale $\hat{f}_{1}(p)$, and functional form $\hat{f}_{2}(p)$ given by,
\begin{subequations}
    \begin{align}
        \hat{f}_{1}(p) &= 
        \begin{cases}
            2p, &0\leq p<0.5 \\
            1+2k_{\gamma}\eta\left( 2p-1 \right), & p \geq 0.5 
        \end{cases} \\
        \hat{f}_{2}(p) &= \begin{cases}
        (2p)^{1/n}, &0\leq p<0.5 \\
        1+2k_{\gamma}\left( 2p-1 \right)^n, & p \geq 0.5 
        \end{cases}
    \end{align}
\end{subequations}
The parameter $\eta$ controls the misspecification in the amplification infection rates due to higher-order interactions, and $n$ affects the degree of nonlinearity from the true scaling function $f(p)$ (see Fig. \ref{fig:scaling_functions}).

\begin{figure}
    \centering
    \includegraphics[width=\linewidth]{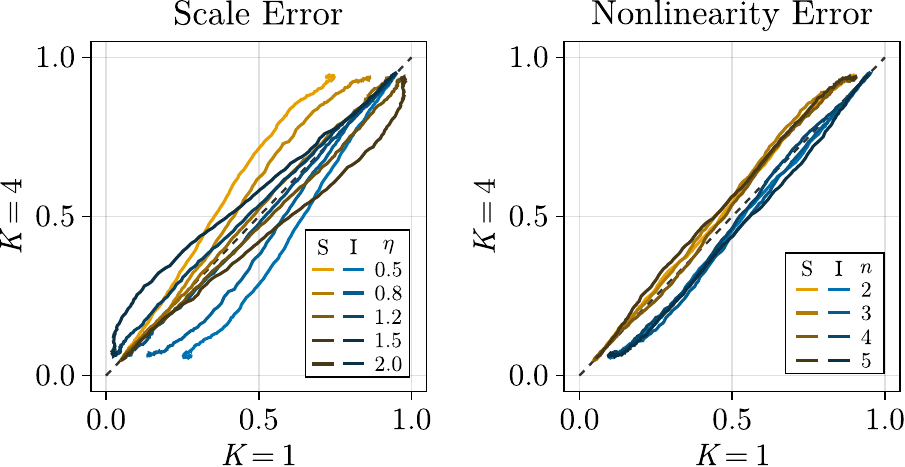}
    \caption{Comparison trajectories showing robust performance across misspecification due to scale and functional form, with the former having more detrimental effects on performance.}
    \label{fig:model_misspec}
\end{figure}

We simulate disease dynamics for trajectory pairs of $K = 1$ and $K=4$ on a 500 node ER random network for $\eta\in\{0.5,0.8,1.2,1.5,2\}$ and $n \in \{2,3,4,5\}$. In both cases, dynamic normalisation of infection rates based on network activity ratios calculated with  $\hat{f}_{1},\hat{f}_{2}$ remained sufficient for pairwise interactions to replicate trajectories from simulations that include higher-order interactions (see Fig. \ref{fig:model_misspec}). We observe that scale errors are more detrimental to the accuracy of the normalisation errors. 

\subsection{Topological Effects}
\subsubsection{Network Topology Assumptions}
All of the cases analysed so far make two key assumptions pertaining to the network topology on which disease dynamics occur: (1) All hypergraphs are defined as clique complexes, (2) macro-scale topological features are uniformly distributed across the network. We first elaborate further on the importance and meaning of these two characteristics, and present additional tests on their impact on the pairwise to higher-order relationship.

The assumption that all hypergraphs belong to a family of clique complexes of some maximal order $K$ places constraints on the existence of hyperedges in the network. Namely, the existence of a given hyperedge $e$ of order $K$ between a set $v(e)=\{ v_i \}_{i=1}^{K+1}$ of $K+1$ nodes in the network implies that all hyperedges of order $k<K$ consisting of a subset of nodes from $v(e)$ also exist. If this condition is met, interactions on hyperedges present the same connectivity as the pairwise case with the only difference being an acceleration of the disease dynamics. While computationally convenient, such an assumption is not necessarily true in real social networks. As an illustrative example, one can consider a sports team consisting of $K+1$ players who regularly engage in close contact as a group by virtue of their vocation. However, sub-groups of players may never interact outside this group setting.

For the second assumption, We clarify ``macro-scale topological features'' to include the various centrality measures of a given node. and ``uniformly distributed'' to mean that the occurrence of a locally significant topological feature (e.g. hubs, high degree) can occur on any node in the network. Therefore, this assumption describes that a given node's connectivity is not conditioned on their neighbours' local topology.

This assumption can be described more generally by considering the deviation local properties from their neighbourhood. For networks, let $\mu(x_i)$ be some measure (in the statistical sense) of a node's property for the $i^{th}$ node. The deviation in measure $\mu$ between a node and it's neighbours is given by,
$$
\Delta \mu(x_i) \;=\; \mu(x_i) - \frac{1}{|\mathcal{N}_{m}(i)|} \sum_{j \in \mathcal{N}_{m}(i)} \mu(x_j),
$$
where $\mathcal{N}_m$ is the set of neighbours for node $i$ with a distance less than or equal to $m$. Large values of $\Delta \mu(x_i)$ implies that the property of a node differs substantially from that of its neighbourhood.

For example, in the case of small-world networks -- such as those generated using the Watts-Strogatz algorithm \cite{watts1998collective} -- shortcut nodes may arise anywhere in the network. For a random Erd\"{o}s-R\'{e}nyi networks, local heterogeneity -- nodes with abnormally high degrees -- have equal probability to occur anywhere within the network. In terms of node deviations, this implies
$$
P\big(\Delta \mu(x_i) = \Delta \mu^*\big),
$$
where $\Delta\mu^*$ is some constant quantity that is independent of $i$.

One can also consider large deviations in terms of distributions. For networks with small deviations, given a local network measure of interest $\mu(x_i)$, and embedding of the network into a space $\{x_i\}_{i=1}^{N} \to \mathcal{M} \subseteq \mathbb{R}^n$, the probability of observing $P(\mu(x_i) = \mu^*)$ in a randomly generated network is uniformly distributed in $\mathcal{M}$.

Examples that do not possess this characteristic are those that arise from preferential-attachment mechanisms such as the Barab\'{a}si-Albert (BA) scale-free network model \cite{barabasi1999emergence}. In this case, high degree nodes (i.e. hubs) will always remain close together in the network due to the preferential attachment mechanism, and thus result in a highly localised distribution of $P(\mu)$ in $\mathcal{M}$.

\subsubsection{Pathological Network Cases}
One point of interest is to test the robustness of the network activity normalisation methods on network topologies to that do not adhere well to the abovementioned topological assumptions. For this, we consider both artificial randomly generated networks, and two different empirical networks taken from contact data that exhibit non-trivial local topologies. For the artificial case, we choose to study a randomly generated BA scale-free $N=300$ node network with additional edges added using the standard ER random network algorithm with varying probability $p \in \{ 0,0.005,0.02 \}$, where $p=0$ corresponds to the standard BA scale-free model of preferential attachment. Increasing values of $p$ dilute the scale-free property of the network resulting in a topology that better adheres to Assumption (2). 

\begin{figure}
\centering
\begin{subfigure}{0.45\textwidth}
    \includegraphics[width=\textwidth]{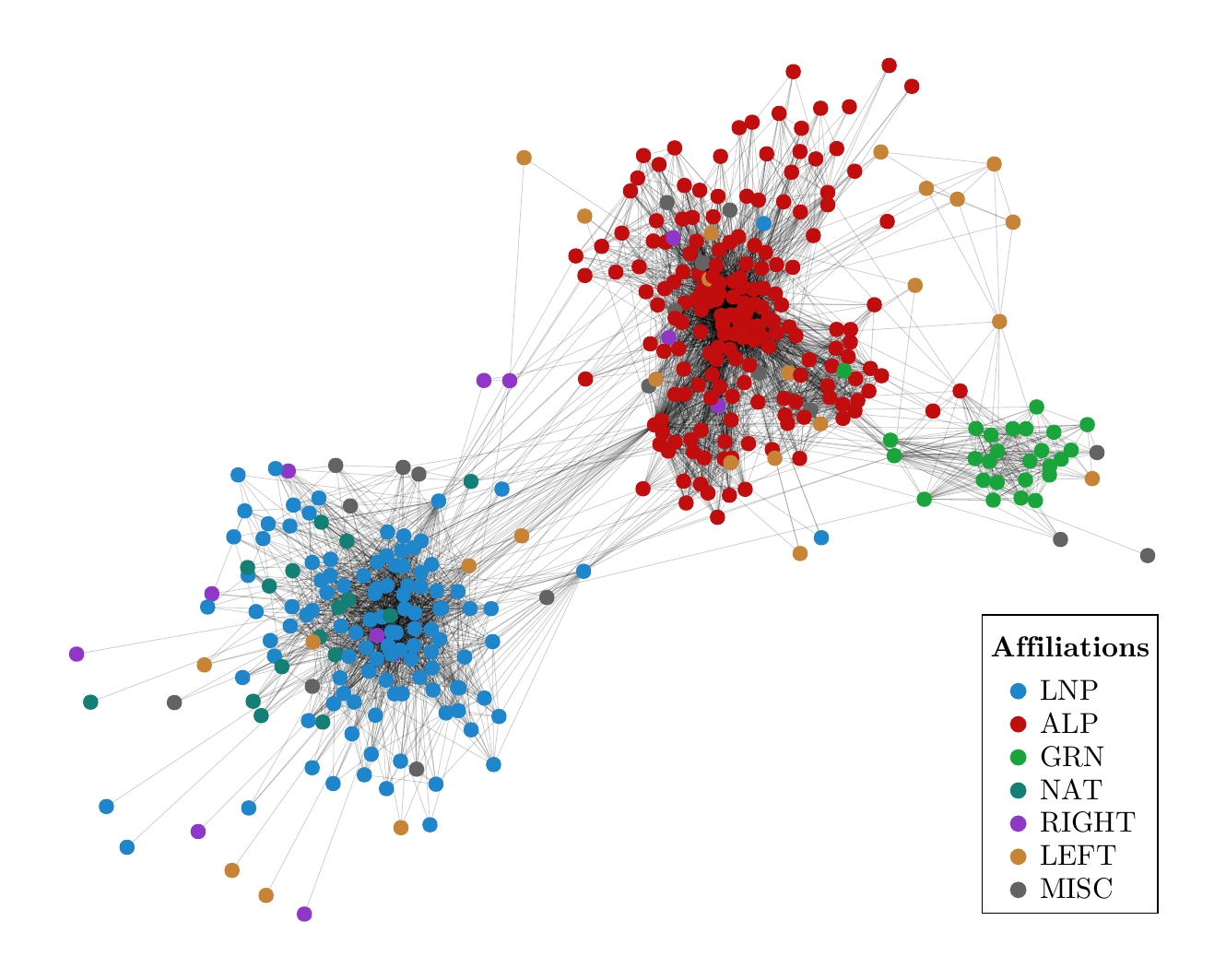}
    \caption{Australian politicians Facebook interactions.}
\end{subfigure}
\hfill
\begin{subfigure}{0.45\textwidth}
    \includegraphics[width=\textwidth]{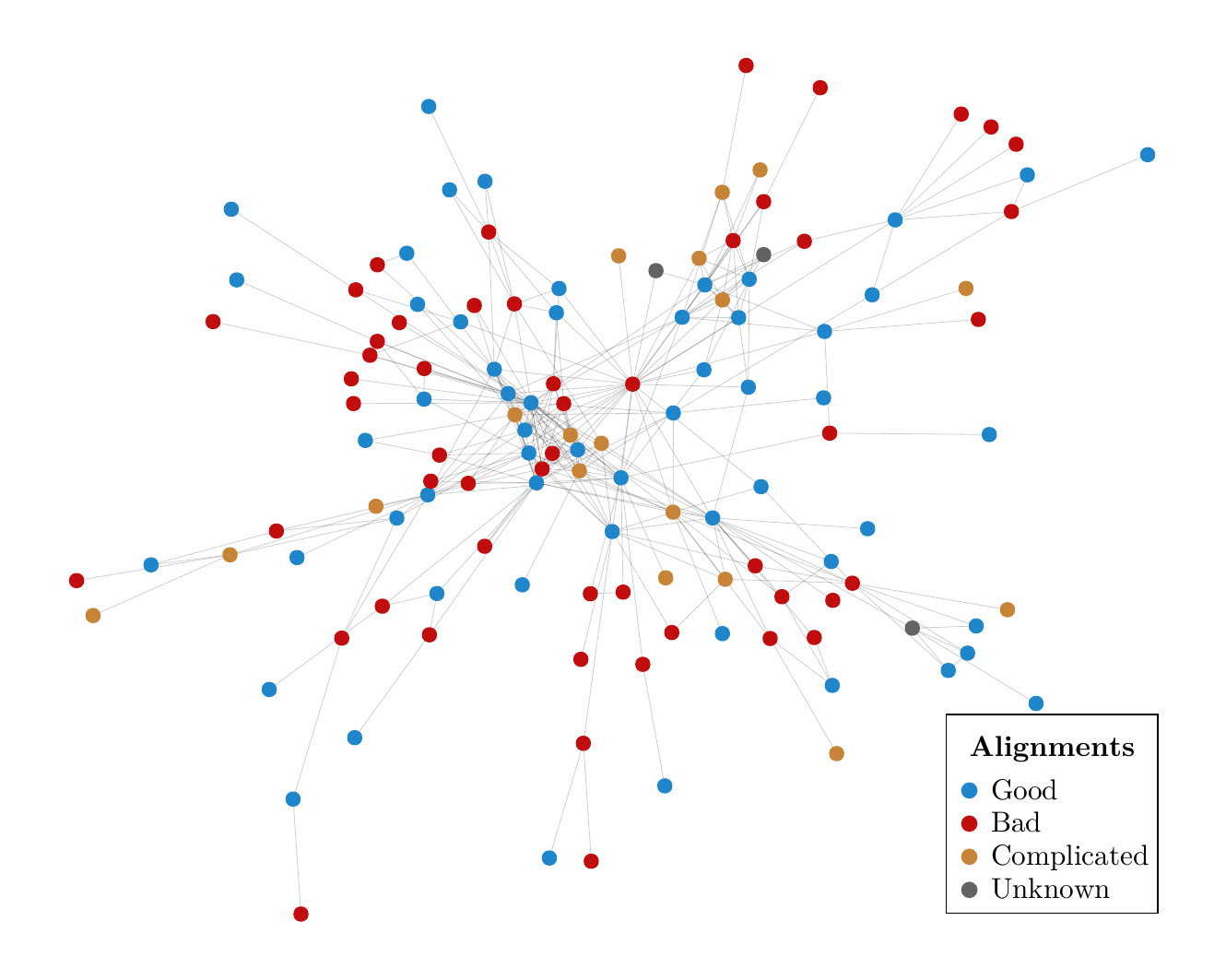}
    \caption{Marvel Cinematic Universe (MCU) conflicts.}
\end{subfigure}
\caption{Visualisation of tested empirical networks with associated node classifications shown. Politician affiliations correspond to either Australian political parties or general political alignment. Hero alignment correspond to the classification of a MCU hero as either a protagonist or antagonist.}
\label{fig:network_plot}
\end{figure}

To demonstrate application to real data, two empirical networks describing human contact were taken from existing publicly available resources. The first is a network of Australian politicians with edges representing mutually liked Facebook pages from November 2017 \cite{NetworkRepository}. The second is a fictional contact network derived from conflicts present from the entire collection of movies from the Marvel Cinematic Universe (MCU) from the Aleph Zero Heroes dataset compiled by Roughan et al. \cite{roughan2020avengers}. Visualisations of both networks are given in Fig. \ref{fig:network_plot}. Both cases showcase features suggestive of preferential attachment with a high degree centre and tree-like or low degree exteriors. 

In addition to examining the impact of base network structures, we also test the importance of the clique complex Assumption (1) by randomly removing hyperedges such that higher order hyperedges no longer contain all their subcliques and thus breaking any implied downward closure of higher order interactions \cite{torres2021and}. This is done using the following pruning algorithm:
\begin{enumerate}
    \item Given an undirected pairwise adjacency matrix $A_0$, construct its clique complex hypergraph $G_K$ up to order $K$.
    \item \label{step:starting_order}Let $k=K$ be the starting maximal order. For each $k^{th}$ order hyperedge $e_k$, delete $e_k$ with probablity $p_d$.
    \begin{enumerate}
        \item If no deletion occurred, continue to the next hyperedge.
        \item If deletion occured, identify all child clique complexes of order $k-1> q > 1$ and recursively repeat for all subsequent children.
        \item If child is a hyperedge of order $q=1$, then only allow the deletion of  theedge if it pairwise network remains connected.
    \end{enumerate}
    \item Repeat step \ref{step:starting_order} with decreasing starting order $k \leftarrow k-1$ until $k = 1$.
\end{enumerate}

An important feature of the above edge deletion algorithm is that lower-order hyperedges experience more deletion events due to the recursive structure of the algorithm. For low to moderate deletion probabilities $p_d$, this results in a skewing of the proportion of hyperedges toward those of higher order.

\subsubsection{Topological Robustness}
\begin{figure*}
    \centering
    \includegraphics[width=\linewidth]{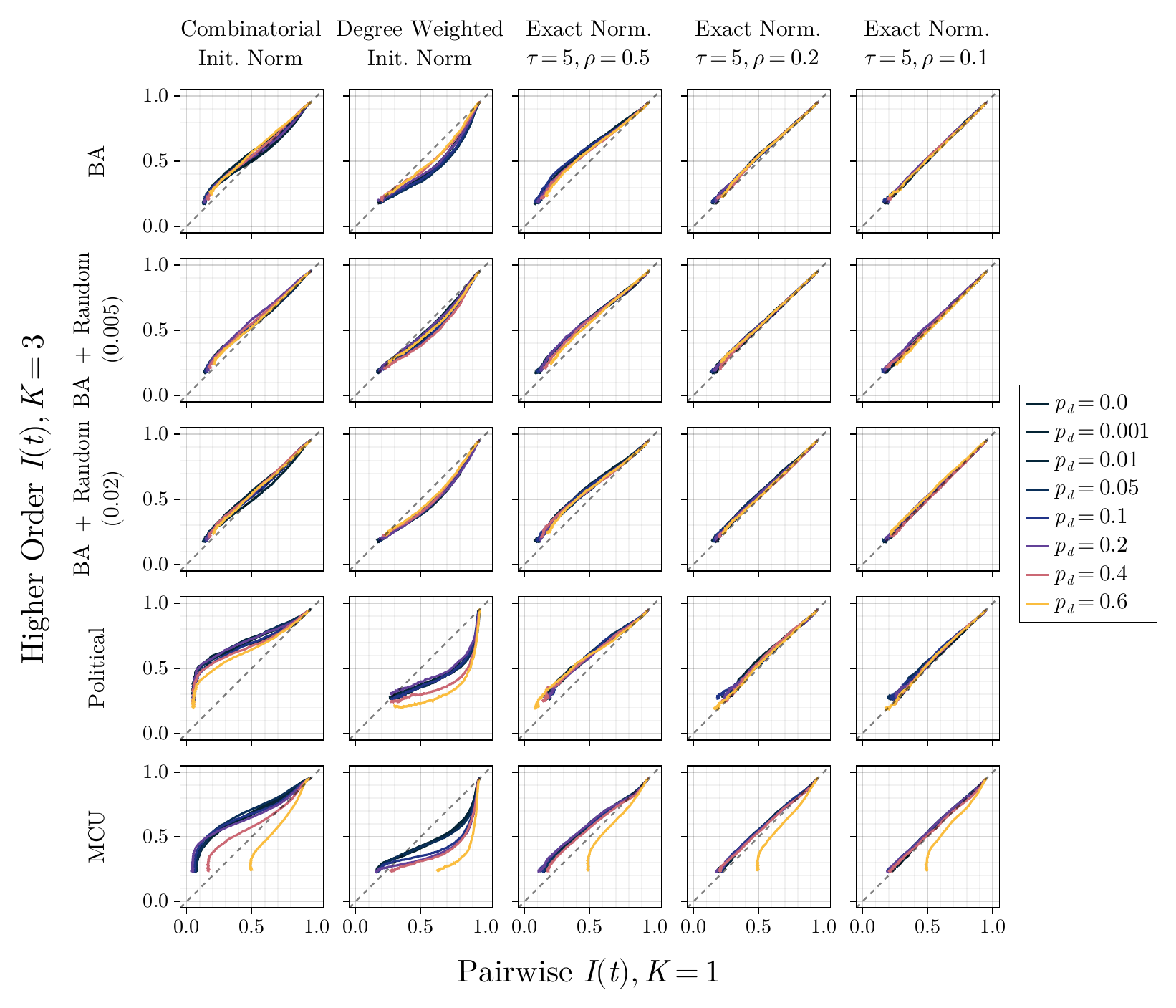}
    \caption{Comparison of pairwise vs higher-order simulation trajectories for 5 tested networks with varying edge deletion probabilities $p_d$.}
    \label{fig:p_deletion}
\end{figure*}

We repeat similar tests as those outlined in Section \ref{subsec:init_norm} and \ref{subsec:adaptive_norm} across the five previously mentioned networks, and $p_d \in (0,0.6)$. Five different normalisation settings are tested consisting of the standard initial normalisation using the combinatorial network activity $\hat{\lambda}$, adaptive normalisation of the exact network activity with decreasing thresholds $\rho \in \{0.5,0.2,0.1 \}$, and a related variant of the combinatorial network activity weighted according to the hyperdegree of network nodes. Results are shown in Fig. \ref{fig:p_deletion}.

In all tested networks, we find that a single initial normalisation of disease parameters using the combinatorial network activity $\hat{\lambda}$ is negatively affected by the highly localised ``macro-scale'' structures all networks such as the clustered groups of hubs. Furthermore, as the influence of higher-order interactions increase (i.e higher $p_d$), the validity of this normalisation also fails. 

To better understand the cause for this failure, a variant of the combinatorial network activity is used as an alternative normalisation. In this case, a surrogate network state is constructed where half of all nodes are assigned an infected state $I$ with priority ranking based on their hyperdegree. The network activity is then calculated for the surrogate network state to yield an approximate network activity $\lambda^* (p_I)$ which subsequently used for normalisation. In all cases, this method was found to correspond to better agreement in the steady state values of the epidemic across all networks.

\begin{figure}
    \centering
    \includegraphics[width=\linewidth]{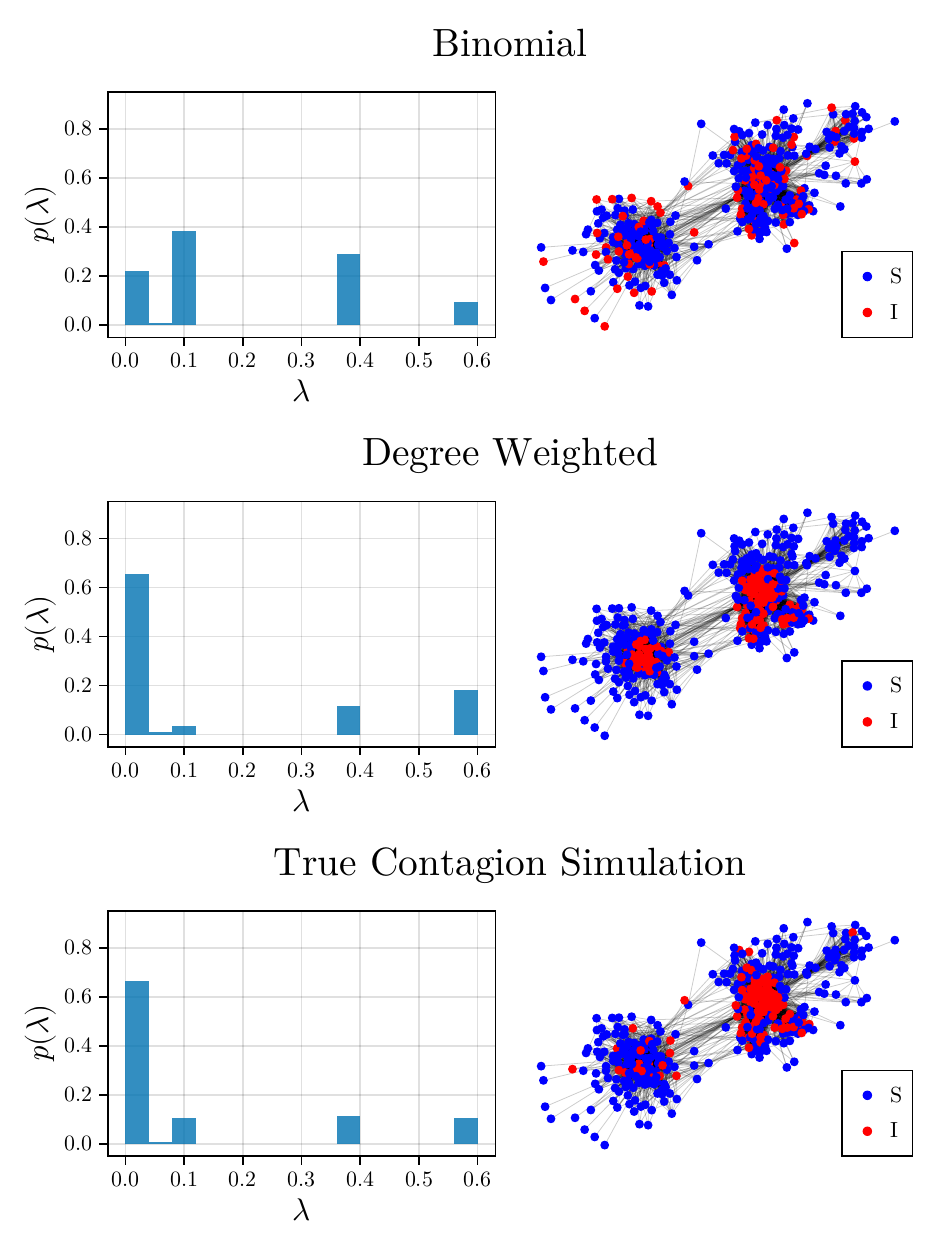}
    \caption{Frequency distributions of network edge activities assuming an infected proportion of $p_I = 0.3$ for the combinatorial (Binomial), degree weighted and true contagion simulation. Surrogate and exact network states are shown on the right. A degree weighted network state surrogate results in a distribution of network activity that more closely resembles the true network activity.}
    \label{fig:network_activity}
\end{figure}

The superior performance of the hyperdegree weighted network activity is due to a closer resemblance between the surrogate network state and the likely true network state, where higher degree nodes will tend to have a higher probability of being infected. This can be seen by looking at the distribution of activities associated with each hyperedge in the combinatorial, degree weighted and true exact network activity cases at a given infection proportion $p_I$ (see Fig. \ref{fig:network_activity}). Notably, the combinatorial approach which inherently assumes a binomial distribution ($p=0.5$) for infections and results in infected nodes that are uniformly spread throughout the entire network. However, node states from contagion dynamics have a tendency to bias infections towards nodes with high degree, which is well replicated by the degree weighting approach.

Interestingly, we note that the degree weighted approach still fails to properly replicate epidemic transients. One may consider a na\"{i}ve extension of the degree weighted network activity that is normalised based on the true current proportions of infected $p_I$ in both the pairwise and higher order cases. This method is attractive as it potentially allows for an approximation of the true exact network activity without requiring knowledge of the full state. However, we find that such an approach fails due to the presence of highly local macro-scale structures that exacerbate the deviation $|\lambda^{*}(p_I) - \lambda|$ as $p_I \to 1$. This deviation can be partially attributed to random recovery of individual nodes within the high-degree clusters of the network, whose hyperedges boost the underlying network activity.

As seen in Fig. \ref{fig:p_deletion}, the exact activity normalisation methods performed well across all tested networks and for moderately high levels of edge deletion, with the exception of $p_d = 0.6$ in the case of the MCU network. This performance is improved for increasingly stringent thresholds $\rho$. However, we note that the performance of the exact method is lower than when applied to the triangular and random networks tested in Section \ref{subsec:adaptive_norm}. This suggests that whilst topologically robust, network activity has a non-trivial relationship with network topology.

\section{\label{sec:level1_conclusion}Conclusion}
In this work, we quantify the effect of higher-order topologies on the spreading dynamics on networks. Focusing on the single disease case, we have presented an agent-based model that unifies the classical compartmental SIR model of homogeneous mixing with pairwise simplicial contagion models on networks, and extends it to account interactions on hyperedges of arbitrary order. We use this model in conjunction with a normalisation method based on the notion of network activity to compare spread dynamics across hypergraphs of different maximal orders.

We first consider the restricted case of clique complexes with uniformly distributed macro-scale topological features for both SIS and SIR dynamics. The inclusion of higher-order topology is found to primarily affect the transient dynamics in epidemic trajectories. Using carefully chosen disease parameters such that initial combinatorial network activities are approximately equal, simulations restricted to pairwise (order 1) interactions are able to replicate the steady state behaviours resulting from higher-order interactions. Furthermore, a sufficiently accurate normalisation can be calculated purely based on the network topology, network state and the scaling function describing the higher-order interaction. We find that allowing disease parameters to dynamically vary over time is sufficient for pairwise simulations to reproduce full trajectories of higher-order interactions. Thus we conclude that there is a duality between topological effects, and the stationarity of disease parameters.

The proposed normalisation method based on network activities makes several assumptions such as having knowledge of the scaling function, and having a clique complex structure with relatively mild macro-scale heterogeneities. While the findings are theoretically interesting, these assumptions impose heavy constraints on applicability as real world systems rarely fulfill these conditions. To this end, we test the robustness of the network activity approach towards model misspecification and non-clique hypergraphs with non-trivial topology taken from empirical data. Overall, all results are relatively robust to model misspecification due to scale and functional form errors. In contrast, we find that network topology plays a larger role where networks with macro-scale heterogeneities -- such as social networks or those arising from preferential attachment mechanisms -- where higher-order interactions is not as easily or robustly approximated by pairwise interactions. This result is particular true for cases where disease parameters are kept static. 

The ability for pairwise interaction models to reproduce epidemic trajectories from models containing higher-order complexity features raises interesting questions on the need for higher-order features to be included in disease models. Furthermore, it encourages a reconsideration of whether dynamics, observed or modelled, should be attributed as the result of complex topological structure or temporally varying spreading mechanisms. The above tests focus on epidemic trajectories, which represent summaries of the disease averaged across the entire network. These results invite further questions on the similarity of spreading dynamics on the local scale of a given network and if the order of infection events between individuals differ between interactions of different order. One can also consider the effect of incubation and immunity refractory periods, which may impart a form of time delay in the infection process resulting in more dynamic disease parameters. On application, network activities offer a potential way to construct adaptive rewiring strategies that can be used to maintain a sustainable level of infection.



\appendix

\bibliography{apssamp}

\end{document}